\documentclass[a4paper,10pt]{amsproc}
\usepackage{amsfonts, amsmath, amssymb, graphicx, mathrsfs, tikz, xcolor}
\usetikzlibrary{angles,quotes}
\tikzset{Bullet/.style={fill=black,draw,color=#1,circle,minimum size=2pt,scale=0.45}}
\usepackage{scalerel,stackengine}
\stackMath
\newcommand\reallywidehat[1]{%
	\savestack{\tmpbox}{\stretchto{%
			\scaleto{%
				\scalerel*[\widthof{\ensuremath{#1}}]{\kern.1pt\mathchar"0362\kern.1pt}%
				{\rule{0ex}{\textheight}}
			}{\textheight}%
		}{2.4ex}}%
	\stackon[-6.9pt]{#1}{\tmpbox}%
}
\usepackage{hyperref}
\hypersetup{
	colorlinks=true,
	linkcolor=blue,
	filecolor=blue,
	citecolor = green,      
	urlcolor=cyan,
}
\usepackage[all]{xy}

\definecolor{lime}{HTML}{A6CE39}
\DeclareRobustCommand{\orcidicon}{
	\begin{tikzpicture}
		\draw[lime, fill=lime] (0,0) 
		circle [radius=0.16] 
		node[white] {{\fontfamily{qag}\selectfont \tiny ID}};
		\draw[white, fill=white] (-0.0625,0.095) 
		circle [radius=0.007];
	\end{tikzpicture}
	\hspace{-2mm}
}

\foreach \x in {A, ..., Z}{\expandafter\xdef\csname orcid\x\endcsname{\noexpand\href{https://orcid.org/\csname orcidauthor\x\endcsname}
		{\noexpand\orcidicon}}
}

\theoremstyle{plain}

\newtheorem{theorem}{Theorem}[section]
\newtheorem{lemma}[theorem]{Lemma}

\theoremstyle{definition}
\newtheorem{definition}[theorem]{Definition}

\newtheorem{remark}[theorem]{Remark}
\newtheorem{counter example}[theorem]{Counter Example}

\newtheorem{corollary}[theorem]{Corollary}
\newtheorem{example}[theorem]{Example}
\newtheorem{convention}[theorem]{Convention}
\numberwithin{equation}{section}

\DeclareMathAlphabet{\mathscr}{OT1}{pzc}{m}{it} 
\begin{document}
	\Large{
		\title{Zero-divisor graph of the rings $C_\mathscr{P}(X)$ and $C^\mathscr{P}_\infty(X)$}
		
		\author[S.K. Acharyya]{Sudip Kumar Acharyya}
		\address{Department of Pure Mathematics, University of Calcutta, 35, Ballygunge Circular Road, Kolkata 700019, West Bengal, India}
		\email{sdpacharyya@gmail.com}

		\author[A. Deb Ray]{Atasi Deb Ray\orcidB}
		\address{Department of Pure Mathematics, University of Calcutta, 35, Ballygunge Circular Road, Kolkata 700019, West Bengal, India}
		\email{debrayatasi@gmail.com}
		
		\author[P. Nandi]{Pratip Nandi\orcidC}
		\address{Department of Pure Mathematics, University of Calcutta, 35, Ballygunge Circular Road, Kolkata 700019, West Bengal, India}
		\email{pratipnandi10@gmail.com}
		
		\thanks{The third author thanks the CSIR, New Delhi – 110001, India, for financial support}
		
\begin{abstract}
	In this article we introduce the zero-divisor graphs $\Gamma_\mathscr{P}(X)$ and $\Gamma^\mathscr{P}_\infty(X)$ of the two rings $C_\mathscr{P}(X)$ and $C^\mathscr{P}_\infty(X)$; here $\mathscr{P}$ is an ideal of closed sets in $X$ and $C_\mathscr{P}(X)$ is the aggregate of those functions in $C(X)$, whose support lie on $\mathscr{P}$. $C^\mathscr{P}_\infty(X)$ is the $\mathscr{P}$ analogue of the ring $C_\infty (X)$. We find out conditions on the topology on $X$, under-which $\Gamma_\mathscr{P}(X)$ (respectively, $\Gamma^\mathscr{P}_\infty(X)$) becomes triangulated/ hypertriangulated. We realize that $\Gamma_\mathscr{P}(X)$ (respectively, $\Gamma^\mathscr{P}_\infty(X)$) is a complemented graph if and only if the space of minimal prime ideals in $C_\mathscr{P}(X)$ (respectively $\Gamma^\mathscr{P}_\infty(X)$) is compact. This places a special case of this result with the choice $\mathscr{P}\equiv$ the ideals of closed sets in $X$, obtained by Azarpanah and Motamedi in \cite{Azarpanah} on a wider setting. We also give an example of a non-locally finite graph having finite chromatic number. Finally it is established with some special choices of the ideals $\mathscr{P}$ and $\mathscr{Q}$ on $X$ and $Y$ respectively that the rings $C_\mathscr{P}(X)$ and $C_\mathscr{Q}(Y)$ are isomorphic if and only if $\Gamma_\mathscr{P}(X)$ and $\Gamma_\mathscr{Q}(Y)$ are isomorphic.
\end{abstract}


\thanks{}
\maketitle
\section{Introduction}
In what follows $X$ stands for a Tychonoff space. Let $\mathscr{P}$ be an ideal of closed sets in $X$ in the following sense: if $A\in\mathscr{P}$ and $B\in\mathscr{P}$, then $A\cup B\in\mathscr{P}$ and if $A\in\mathscr{P}$ and $C\subset A$ with $C$, closed in $X$, then $C\in\mathscr{P}$. Suppose $C_\mathscr{P}(X)$ is the family of all those function $f$ in $C(X)$ whose support $cl_X(X\setminus Z(f))\in\mathscr{P}$, here $Z(f)=\{x\in X:f(x)=0\}$ is the zero set of $f$. Suppose $C^\mathscr{P}_\infty(X)=\{f\in C(X):\text{for each }\epsilon>0,\{x\in X:|f(x)|\geq\epsilon\}\in\mathscr{P}\}$. It turns out that $C_\mathscr{P}(X)$ and $C^\mathscr{P}_\infty(X)$ are both commutative rings, possibly without identity and $C_\mathscr{P}(X)\subset C^\mathscr{P}_\infty(X)$. Let $\Gamma(C_\mathscr{P}(X))\equiv \Gamma_\mathscr{P}(X)$ be the graph, whose vertices are non-zero divisors of zero in $C_\mathscr{P}(X)$ and $\Gamma^\mathscr{P}_\infty(X)$ be the analogous graph associated with $C^\mathscr{P}_\infty(X)$. Two distinct vertices $f$ and $g$ in $C_\mathscr{P}(X)$ (respectively in $C^\mathscr{P}_\infty(X)$) are said to be connected by an edge, in which case they are called adjacent vertices if and only if $f.g=0$. Our intention to write this article is to establish a number of facts which highlight possible interaction between graph properties of $\Gamma_\mathscr{P}(X)$ (respectively $\Gamma^\mathscr{P}_\infty(X)$) and ring properties of $C_\mathscr{P}(X)$ (respectively $C^\mathscr{P}_\infty(X)$) leading to further interaction between these two properties and the topological properties on $X$. It is easy to see that on choosing $\mathscr{P}$ to be the ideal of all closed sets in $X$, $C_\mathscr{P}(X)$ becomes identical to $C(X)$. We realize that some of the results related to zero-divisor graph of $C(X)$ obtained in \cite{Azarpanah} are special cases of facts obtained in the present paper. We would like to mention in this context as far as we dig into literature that there are only four papers on graphs having their vertices lying in $C(X)$. See the articles \cite{Alafifi}, \cite{Azarpanah}, \cite{Badie} and \cite{Bose} in this context. In the technical section \ref{Sec2} of this paper we introduce several well-known parameters related to the graph $\Gamma_\mathscr{P}(X)$. These include distance between distinct vertices, diameter, radius of a graph, eccentricity of a vertex. We show that, the distance between any two vertices of the graph $\Gamma_\mathscr{P}(X)$ is at most $3$. This leads to necessary and sufficient condition for $\Gamma_\mathscr{P}(X)$ to be triangulated (respectively hypertriangulated).\\
In section \ref{Sec3}, we compute the lengths of various possible cycles and determine the eccentricity of the vertices in some case. We also find out a few relations interconnecting the dominating number and clique number of $\Gamma_\mathscr{P}(X)$ and the cellularity of $X$. Furthermore we determine when does $\Gamma_\mathscr{P}(X)$ becomes a complemented graph.\\
In section \ref{Sec4}, we calculate several parameters related to the graph $\Gamma^\mathscr{P}_\infty(X)$. These are mostly parallel to their $\Gamma_\mathscr{P}(X)$ analogues obtained in the section \ref{Sec2} and \ref{Sec3}. However we proved that the chromatic number of $\Gamma_\mathscr{P}(X)$ and $\Gamma^\mathscr{P}_\infty(X)$ are identical.\\
If two rings $C_\mathscr{P}(X)$ and $C_\mathscr{Q}(Y)$ are isomorphic, then it is easy to see that their zero-divisor graphs $\Gamma_\mathscr{P}(X)$ and $\Gamma_\mathscr{Q}(Y)$ are isomorphic in the following sense: there is a bijection between the set of vertices of these two graphs which preserve the adjacency relation. However the converse problem to find out any possible isomorphism between the rings $C_\mathscr{P}(X)$ and $C_\mathscr{Q}(Y)$ on the basis of the hypothesis that there is a graph isomorphism between $\Gamma_\mathscr{P}(X)$ and $\Gamma_\mathscr{Q}(Y)$, in general appears to be too wild to venture into. Nevertheless, by making some special choices of ideals $\mathscr{P}$ and $\mathscr{Q}$ on $X$ and $Y$ respectively, we make some breakthrough in this matter. We establish that if $\mathscr{P}$ is the ideal of all finite subsets of $X$ and $\mathscr{Q}$, the ideal of all finite subsets of $Y$, then the rings $C_\mathscr{P}(X)$ and $C_\mathscr{Q}(Y)$ are isomorphic if and only if $\Gamma_\mathscr{P}(X)$ and $\Gamma_\mathscr{Q}(Y)$ are isomorphic [Theorem \ref{Th6.4}]. This is the final result in section \ref{Sec6}.\\
For more information on the rings $C_\mathscr{P}(X)$ and $C^\mathscr{P}_\infty(X)$, the reader is refereed to see the articles \cite{Acharyya} and \cite{Acharyya1}. For graph theoretic information, the reader is refereed to the book \cite{Diestel}. 
\section{Technical notations related to $\Gamma_\mathscr{P}(X)$}\label{Sec2}
The distance between two distinct vertices $f$ and $g$ in $\Gamma_\mathscr{P}(X)$, denoted by $d(f,g)$, is the length of the shortest path from $f$ to $g$. We wish to denoted by $V_\mathscr{P}(X)$, the set of vertices of the graph $\Gamma_\mathscr{P}(X)$. The diameter of the graph $\Gamma_\mathscr{P}(X)$ is defined by: $diam(\Gamma_\mathscr{P}(X))=Max\{d(f,g):f,g\in V_\mathscr{P}(X)\}$. The eccentricity $e(f)$ of an $f\in V_\mathscr{P}(X)$ is defined by: $e(f)=Max\{d(f,g):g\in V_\mathscr{P}(X)\}$. An $f\in V_\mathscr{P}(X)$ is called a center of $\Gamma_\mathscr{P}(X)$ if $e(f)\leq e(g)$ holds for each $g\in V_\mathscr{P}(X)$ and in this case $e(f)$ is called the radius of the graph. The girth of $\Gamma_\mathscr{P}(X)$, denoted by $gr(\Gamma_\mathscr{P}(X))$, is the length of the smallest cycle in this graph. Like any graph $\Gamma_\mathscr{P}(X)$ is called triangulated (respectively hypertriangulated) if each vertex (respectively each edge) of this graph is a vertex (respectively is an edge) of a triangle. The smallest length of a cycle containing two distinct vertices $f$ and $g$ in $\Gamma_\mathscr{P}(X)$ will be denoted by $c(f,g)$.\\
A subset $D$ of $V_\mathscr{P}(X)$ is called a dominating set in $\Gamma_\mathscr{P}(X)$ if for each $f\in V_\mathscr{P}(X)\setminus D$, there exists $g\in D$ such that $f$ and $g$ are adjacent. The dominating number of $\Gamma_\mathscr{P}(X)$ is defined as follows: $dt(\Gamma_\mathscr{P}(X))=min \{|D|:D\text{ is a dominating set in }\Gamma_ \mathscr{P}(X)\}$. A coloring of a graph is a labeling of the vertices of the graph with colors such that no two adjacent vertices have the same color. More precisely, for a cardinal number $\alpha$ (finite or infinite), an $\alpha$-coloring of $\Gamma_\mathscr{P}(X)$ is a map $\psi:V_\mathscr{P}(X)\to[0,\alpha)$ with the following condition: whenever $f,g\in V_\mathscr{P}(X)$ and $f.g=0$, $\psi(f)\neq\psi(g)$. The chromatic number of $\Gamma_\mathscr{P}(X)$ is defined as follows: $\chi(\Gamma_\mathscr{P}(X)) =min\{\alpha:\text{there exists a }\alpha\text{-coloring of }\Gamma_ \mathscr{P}(X)\}$.\\
A complete subgraph of a graph $G$ is any subset $H$ of $G$ such that each pair of distinct vertices in $H$ are adjacent. The clique number of $G$ is defined as follows: $\omega(G)=sup\{|H|\text{ is a complete subgraph of }G\}$.\\
A collection $\mathscr{B}$ of non-empty open sets in $X$ is called a cellular family if any two distinct members of $\mathscr{B}$ are disjoint. The cellularity of a space $X$ is defined as follows: $c(X)=sup\{|\mathscr{B}|:\mathscr{B}\text{ is a cullular family of open}$  $\text{sets in }X\}$.
\begin{definition}
	$X$ is called locally $\mathscr{P}$ at a point $x\in X$, if there exists an open neighbourhood $V$ of $x$ in $X$ such that $cl_XV\in\mathscr{P}$. $X$ is said to be locally $\mathscr{P}$ if it is locally $\mathscr{P}$ at each point on it.
\end{definition}
Let $X_\mathscr{P}=\{x\in X: X\text{ is  locally }\mathscr{P}\text{ at }x\}$. Then it is easy to prove that $X_\mathscr{P}$ is an open set in $X$. Also,  $X$ is locally $\mathscr{P}$ if and only if $X_\mathscr{P}=X$.
\begin{lemma}\label{Lem2.2}
	Given $x\in X_\mathscr{P}$ and an open neighbourhood $G$ of $x$, there exists $f\in C_\mathscr{P}(X)$ such that $x\in X\setminus Z(f)\subset cl_X(X\setminus Z(f))\subset G$.
\end{lemma}
\begin{proof}
	Since $x\in X_\mathscr{P}$ there exists an open neighbourhood $U$ of $x$ in $X$ such that $cl_XU\in\mathscr{P}$. Consider the open neighbourhood $U\cap G$ of $x$. Then by complete regularity of $X$, there exists $f\in C(X)$ such that $x\in X\setminus Z(f)\subset cl_X(X\setminus Z(f))\subset U\cap G$. Since $cl_X(X\setminus Z(f))\subset U\subset cl_XU\in\mathscr{P}$, it follows that $cl_X(X\setminus Z(f))\in\mathscr{P}$, i.e., $f\in C_\mathscr{P}(X)$.
\end{proof}
The following result decides which non-zero elements in $C_\mathscr{P}(X)$ are vertices in the graph $\Gamma_\mathscr{P}(X)$.
\begin{theorem}\label{Th2.3}
	For any $f\in C_\mathscr{P}(X)\setminus\{0\}$, the following three statements are equivalent:
	\begin{enumerate}
		\item $f\in V_\mathscr{P}(X)$
		\item $X_\mathscr{P}-cl_X(X\setminus Z(f))= X_\mathscr{P}\cap int_X(Z(f))\neq\emptyset$
		\item $cl_X(X_\mathscr{P})\cap int_XZ(f)\neq\emptyset$
	\end{enumerate}
\end{theorem}
\begin{proof}
	\underline{$(1)\implies(2)$:} Let $(1)$ hold. Then there exists $g\in V_\mathscr{P}(X)$ such that $f.g=0$. This implies that $X\setminus Z(f)\cap X\setminus Z(g)=\emptyset$. Choose a point $x\in X\setminus Z(g)$. Then $x\in X_\mathscr{P}$ because $X\setminus Z(g)\subset X_\mathscr{P}$. On the other hand $x\notin cl_X(X\setminus Z(f))$, i.e., $x\in int_XZ(f)$. Thus $x\in X_\mathscr{P}\cap int_XZ(f)$.\\
	\underline{$(2)\implies(3)$:} This is trivial.\\
	\underline{$(3)\implies(1)$:} Let $(3)$ be true. Choose a point $p\in cl_X(X_\mathscr{P})-cl_X(X\setminus Z(f))$ i.e., $p\in cl_X(X_\mathscr{P})\cap int_XZ(f)$. Clearly then $X_\mathscr{P}\cap int_XZ(f)\neq\emptyset$. Choose a point $q\in X_\mathscr{P}\cap int_XZ(f)$. Then by Lemma \ref{Lem2.2}, there exists $g\in C_\mathscr{P}(X)$ such that $q\in X\setminus Z(g)\subset cl_X(X\setminus Z(g))\subset int_XZ(f)$. Consequently $f.g=0$ and $g\in C_\mathscr{P}(X)$ is a vertex in $\Gamma_\mathscr{P}(X)$. Thus $f\in V_\mathscr{P}(X)$.
\end{proof}
\begin{corollary}\label{Cor2.4}
	Let $cl_X(X_\mathscr{P})\notin\mathscr{P}$. Then each non-zero element $f$ of $C_\mathscr{P}(X)$ is a vertex of $\Gamma_\mathscr{P}(X)$.
\end{corollary}
\begin{proof}
	If possible let there exist $f\in C_\mathscr{P}$, $f\neq 0$ such that $f\notin V_\mathscr{P}(X)$. Then $cl_X(X_\mathscr{P})\cap int_XZ(f)=\emptyset$ and consequently $cl_X(X_\mathscr{P})\subset X\setminus int_XZ(f)=cl_X(X\setminus Z(f))$. Since $cl_X(X\setminus Z(f))\in\mathscr{P}$, this implies that $cl_X(X_\mathscr{P})\in\mathscr{P}$, a contradiction.
\end{proof}
The converse of the last corollary is not true. The following is a simple counterexample.
\begin{example}
	Let $X=\mathbb{Q}\equiv$ the space of all rational numbers. Suppose $\mathscr{P}\equiv$ the ideal of all compact subsets of $X$. Since $\mathbb{Q}$ is nowhere locally compact, it follows that $X_\mathscr{P}=\emptyset$. Consequently, $cl_XX_\mathscr{P}=\emptyset\in\mathscr{P}$. But since $C_\mathscr{P}(X)=C_K(\mathbb{Q})=\{0\}$ [see $4D2$, \cite{GJ}]. The condition that each non-zero element in $C_\mathscr{P}(X)$ is a member of $\Gamma_\mathscr{P}(X)$ is vacuously satisfied. 
\end{example}
We shall show in Section \ref{Sec6} that by an appropriate choice of the space $X$ and the ideal $\mathscr{P}$ of closed sets in $X$, the converse of Corollary \ref{Cor2.4} is true.\\
The following proposition decides, when a given pair of vertices in $\Gamma_\mathscr{P}(X)$ admit of a third vertex adjacent to both of them.
\begin{theorem}\label{Th2.6}
	Let $f,g\in V_\mathscr{P}(X)$. Then there exists a vertex $h$, adjacent to both $f$ and $g$ if and only if $int_X(Z(f)\cap Z(g)\cap X_\mathscr{P})\neq\emptyset$.
\end{theorem}
\begin{proof}
	First let $h\in V_\mathscr{P}(X)$ be adjacent to $f$ and $g$. Then $h\in C_\mathscr{P}(X)$ implies that $X\setminus Z(h)\subset X_\mathscr{P}$. On the other hand $h.f=0=g.h$ implies that $\emptyset\neq X\setminus Z(h)\subset Z(f)\cap Z(g)$. Thus $\emptyset\neq X\setminus Z(h)\subset X_\mathscr{P}\cap Z(f)\cap Z(g)$. Hence $int_X(Z(f)\cap Z(g)\cap X_\mathscr{P})\neq\emptyset$.\\
	Conversely let $int_X(Z(f)\cap Z(g)\cap X_\mathscr{P})\neq\emptyset$. So there exists a non-empty open set $W$ contained in $Z(f)\cap Z(g)\cap X_\mathscr{P}$. Choose a point $x\in W$. Then by Lemma \ref{Lem2.2}, there exists $h\in C_\mathscr{P}(X)$ such that $x\in X\setminus Z(h)\subset cl_X(X\setminus Z(h))\subset W\subset Z(f)\cap Z(g)\cap X_\mathscr{P}$. It is clear that $h\neq 0$ and $h.f=h.g=0$. Then $h\in V_\mathscr{P}(X)$ and is adjacent to both $f$ and $g$.
\end{proof}
In the next theorem, we compute the possible distance between pairs of distinct vertices.
\begin{theorem}\label{Th2.7}
	Let $f,g\in V_\mathscr{P}(X)$. Then 
	\begin{enumerate}
		\item $d(f,g)=1$ if and only if $X\setminus Z(f)\cap X\setminus Z(g)=\emptyset$.\label{Th2.7(1)}
		\item $d(f,g)=2$ if and only if $X\setminus Z(f)\cap X\setminus Z(g)\neq\emptyset$ and $int_X(Z(f)\cap Z(g)\cap X_\mathscr{P})\neq\emptyset$.\label{Th2.7(2)}
		\item $d(f,g)=3$ if and only if $X\setminus Z(f)\cap X\setminus Z(g)\neq\emptyset$ and $int_X(Z(f)\cap Z(g)\cap X_\mathscr{P})=\emptyset$.\label{Th2.7(3)}
	\end{enumerate}
\end{theorem}
\begin{proof}
	\hspace{3cm}
	\begin{enumerate}
		\item Trivial
		\item Follows immediately from part (\ref{Th2.7(1)}) of this Theorem and Theorem \ref{Th2.6}.
		\item Suppose $d(f,g)=3$. Then it follows from part (\ref{Th2.7(1)}) and (\ref{Th2.7(2)}) of the present Theorem that $X\setminus Z(f)\cap X\setminus Z(g)\neq\emptyset$ and $int_X(Z(f)\cap Z(g)\cap X_\mathscr{P})=\emptyset$.\\
		Conversely let $X\setminus Z(f)\cap X\setminus Z(g)\neq\emptyset$ and $int_X(Z(f)\cap Z(g)\cap X_\mathscr{P})=\emptyset$. Then it follows from part (\ref{Th2.7(1)}) and (\ref{Th2.7(2)}) that $d(f,g)>2$. Since $f,g\in V_\mathscr{P}(X)$, there exist $f_1,g_1\in V_\mathscr{P}(X)$ such that $f.f_1=0=g.g_1$. To ascertain that $d(f,g)=3$, it suffices to show that $f_1.g_1=0$. Indeed $f.f_1=0$ implies that $X\setminus Z(f_1) \subset X_\mathscr{P}-cl_X(X\setminus Z(f))=X_\mathscr{P}\cap int_XZ(f)$. Analogously $X\setminus Z(g_1)\subset X_\mathscr{P}\cap int_XZ(g)$. It follows that: $X\setminus Z(f_1)\cap X\setminus Z(g_1)\subset X_\mathscr{P}\cap int_XZ(f)\cap int_XZ(g)=int_X(X_\mathscr{P}\cap Z(f)\cap Z(g))$. The hypothesis $int_X(Z(f)\cap Z(g)\cap X_\mathscr{P})=\emptyset$, therefore implies that $X\setminus Z(f_1)\cap X\setminus Z(g_1)=\emptyset$. Hence $f_1.g_1=0$.
	\end{enumerate}
\end{proof}
\begin{corollary}\label{Cor2.8}
	If $diam(\Gamma_\mathscr{P}(X))=3$, then $cl_X(X_\mathscr{P})\in\mathscr{P}$.
\end{corollary}
\begin{proof}
	Let $diam(\Gamma_\mathscr{P}(X))=3$. So there exists a pair of vertices $f,g\in V_\mathscr{P}(X)$ such that $d(f,g)=3$. Then from Theorem $2.7(\ref{Th2.7(3)})$ it follows that, $int_X(X_\mathscr{P}\cap Z(f)\cap Z(g))=\emptyset$. Now $X_\mathscr{P}-(cl_X(X\setminus Z(f))\cup cl_X(X\setminus Z(g)))=X_\mathscr{P}\cap int_XZ(f)\cap int_XZ(g)=int_X(X_\mathscr{P}\cap Z(f)\cap Z(g))=\emptyset$, i.e., $X_\mathscr{P}\subset cl_X(X\setminus Z(f))\cup cl_X(X\setminus Z(g))$ implies $cl_X(X_\mathscr{P})\subset cl_X(X\setminus Z(f))\cup cl_X(X\setminus Z(g))\in\mathscr{P}$ and hence $cl_X(X_\mathscr{P})\in\mathscr{P}$.
\end{proof}
We would like to mention in this context, the following result proved in \cite{Azarpanah}, Corollary $1.3$ on choosing $\mathscr{P}\equiv$ the ideal of all closed sets in $X$.
\begin{theorem}
	Whenever $X$ has atleast three points, then the diameter of the zero-divisor graph of $C(X)$ is $3$.
\end{theorem}
Thus converse of the Corollary \ref{Cor2.8} is not true. Consider the following example
\begin{example}
	Let $X$ be a Tychonoff space with $|X|=2$ and $\mathscr{P}\equiv$ the ideal of all closed sets in $X$. Then $C_\mathscr{P}(X)=C(X)$ and $X_\mathscr{P}=X$. It can be easily proved that the zero-divisor graph of $C(X)$ is complete bipartite and hence $diam(\Gamma_\mathscr{P}(X))=2$. But $cl_X(X_\mathscr{P})=X\in\mathscr{P}$.
\end{example}
Before proceeding further we have to rule out the cases where $V_\mathscr{P}(X)=\emptyset$, i.e., where the graph is empty.
\begin{theorem}
	For an ideal $\mathscr{P}$ of closed sets in $X$, the following two statements are equivalent:
	\begin{enumerate}
		\item $V_\mathscr{P}(X)\neq\emptyset$.
		\item $|X_\mathscr{P}|\geq 2$.
	\end{enumerate}
\end{theorem}
\begin{proof}
	If $X_\mathscr{P}=\emptyset$, i.e., $X$ is nowhere locally $\mathscr{P}$, then from Theorem \ref{Th2.3}, it follows that no non-zero element in $C_\mathscr{P}(X)$ can be a vertex, in other words $V_\mathscr{P}(X)=\emptyset$. On the other hand, if $X$ is locally $\mathscr{P}$ just at a single point $p$ on $X$, then $p$ is an isolated point of $X$ and $\{p\}\in\mathscr{P}$. Now if $f$ is a non-zero function in $C_\mathscr{P}(X)$ and $f\in V_\mathscr{P}(X)$, then from Theorem \ref{Th2.3}, it follows that $p\in int_XZ(f)$. Consequently, $cl_X(X\setminus Z(f))(\subset X\setminus\{p\})\in\mathscr{P}$ and any point on $X\setminus Z(f)$ is a member of $X_\mathscr{P}$ -- a contradiction to the initial assumption that $X$ is locally $\mathscr{P}$ only at the point $p$. Thus $f$ can not be a non-zero function in $C_\mathscr{P}(X)$, in other words $f\equiv 0$. Hence $V_\mathscr{P}(X)=\emptyset$. Then $\underline{(1)\implies(2)}$ is proved.\\
	$\underline{(2)\implies(1):}$ Assume that $(2)$ is true. So we can find out a pair of distinct points $p,q$ from $X_\mathscr{P}$. Then there exists a co-zero set neighbourhood $C_p$ of $p$ in $X$ such that $cl_X(C_p)\in\mathscr{P}$ and also there is a co-zero set neighbourhood $C_q$ of $q$ in $X$ with $cl_X(C_q)\in\mathscr{P}$. By using the complete regularity of $X$, we can find out co-zero set neighbourhoods $C_p^*$ and $C_q^*$ of $p$ and $q$ in $X$ respectively such that $C_p^*\cap C_q^*=\emptyset$. Let $\reallywidehat{C_p}=C_p^*\cap C_p$ and $\reallywidehat{C_q}=C_q^*\cap C_q$. Then $\reallywidehat{C_p}$ and $\reallywidehat{C_q}$ are disjoint co-zero set neighbourhoods of $p$ and $q$ respectively with $cl_X(\reallywidehat{C_p})\in\mathscr{P}$ and $cl_X(\reallywidehat{C_q})\in \mathscr{P}$. We can write $\reallywidehat{C_p}=X\setminus Z(f)$ for some $f\in C(X)$. Then $f\in C_\mathscr{P}(X)$ as $cl_X(X\setminus Z(f))\in\mathscr{P}$ and also $f\neq 0$. Furthermore $\reallywidehat{C_q}\subset Z(f)$ with $q\in X_\mathscr{P}$. This shows that $q\in int_XZ(f)$. Thus $X_\mathscr{P}\cap int_XZ(f)\neq\emptyset$. Hence from Theorem \ref{Th2.3}, we get that $f\in V_\mathscr{P}(X)$. Therefore, $V_\mathscr{P}(X)\neq\emptyset$.
\end{proof}
\begin{convention}\label{Con2.9}
	In what follows, we shall assume that $|X_\mathscr{P}|\geq 2$ and this will ensure that the graph $\Gamma_\mathscr{P}(X)$ will be non-void.
\end{convention}
The following result shows that for some choice of $\mathscr{P}$, each vertex of $\Gamma_\mathscr{P}(X)$ will be a center, i.e., $\Gamma_\mathscr{P}(X)$ is a self-centric graph. [A graph is said to be a self-centric graph if every vertex in the graph is a center.]
\begin{theorem}\label{Th2.10}
	Let $cl_X(X_\mathscr{P})\notin\mathscr{P}$. Then for an arbitrary $f\in V_\mathscr{P}(X)$, $e(f)=2$.
\end{theorem}
\begin{proof}
	From Corollary \ref{cor2.4}, we get that each non-zero function in $C_\mathscr{P}(X)$ is a vertex of $\Gamma_\mathscr{P}(X)$. So if $g\in C_\mathscr{P}(X)$ , then $d(f,g)=1$ if $f$ and $g$ are adjacent. Suppose $f$ and $g$ are not adjacent. Then since $f^2+g^2$ is a vertex of $\Gamma_\mathscr{P}(X)$, it follows from Theorem \ref{Th2.3} that $X_\mathscr{P}\cap int_XZ(f^2+g^2)\neq\emptyset$, in other words: $X_\mathscr{P}\cap int_XZ(f)\cap int_XZ(g)\neq\emptyset$, i.e., $int_X(X_\mathscr{P}\cap Z(f)\cap Z(g))\neq\emptyset$, it follows from Theorem \ref{Th2.7} that $d(f,g)=2$. Now it is clear that for $f\in V_\mathscr{P}(X)$, $f$ and $2f$ are not adjacent, hence $d(f,2f)=2$. Thus $e(f)=2$.
\end{proof}
\begin{corollary}
	For a locally compact non-compact space $X$, each vertex of the zero-divisor graph of the ring $C_K(X)$ of all continuous functions with compact support is a center of it and the radius of the graph is $2$.
\end{corollary}
\begin{remark}
	We show that the conclusion of Theorem \ref{Th2.10} may not be valid if the hypothesis $cl_X(X_\mathscr{P})\notin\mathscr{P}$ is dropped. Let $X=\mathbb{R}$ and $\mathscr{P}\equiv$ the ideal of closed sets in $\mathbb{R}$. Then $X_\mathscr{P}=X=\mathbb{R}$ and therefore $cl_X(X_\mathscr{P})=\mathbb{R}\in\mathscr{P}$. In this case $C_\mathscr{P}(X)=C(X)=C(\mathbb{R})$. Let $f$ be any vertex of the zero-divisor graph of $C(\mathbb{R})$. Then $int_\mathbb{R}Z(f)\neq\emptyset$. Choose a point $x\in\mathbb{R}\setminus Z(f)$. Then there exists a zero set neighbourhood $Z(g)$ of $x$ in $\mathbb{R}$ such that $Z(g)\cap Z(f)=\emptyset$. We see that $g$ is a vertex of this graph and $int_\mathbb{R}(X_\mathscr{P}\cap Z(f)\cap Z(g))=int_\mathbb{R}(Z(f)\cap Z(g))=\emptyset$. Furthermore $X\setminus Z(f)\cap X\setminus Z(g)\neq\emptyset$. It follows from Theorem $2.7(\ref{Th2.7(3)})$ that $d(f,g)=3$. Then $e(f)=3$. We note that each vertex in this graph is a center but its radius is $3$.
\end{remark}
The following result is a key one to characterize triangulated graph of the form $\Gamma_\mathscr{P}(X)$.
\begin{theorem}\label{Th2.13}
	Let $f\in V_\mathscr{P}(X)$. Then $f$ is a vertex of a triangle if and only if $|X_\mathscr{P}\cap int_XZ(f)|\geq 2$.
\end{theorem}
\begin{proof}
	First let $f$ be a vertex of a triangle. Then there exist $g,h\in V_\mathscr{P}(X)$ such that $f.g=g.h=h.f=0$. Choose $x\in X\setminus Z(g)$ and $y\in X\setminus Z(h)$. Then $\{x,y\}\subset int_XZ(f)\cap X_\mathscr{P}$. Since $X\setminus Z(g)\cap X\setminus Z(h)=\emptyset$, it follows that $x\neq y$. Thus $|X_\mathscr{P}\cap int_XZ(f)|\geq 2$.\\
	Conversely let $|X_\mathscr{P}\cap int_XZ(f)|\geq 2$. Choose $x,y\in X_\mathscr{P}\cap int_XZ(f)$, $x\neq y$. Then $X_\mathscr{P}-( cl_X(X\setminus Z(f))\cup\{y\})$ is an open neighbourhood of $x$ in $X$. By Lemma \ref{Lem2.2}, there exists $g\in V_\mathscr{P}(X)$ such that $x\in X\setminus Z(g)\subset cl_X(X\setminus Z(g))\subset X_\mathscr{P}-(cl_X(X\setminus Z(f))\cup\{y\})$. Clearly $g.f=0$. Furthermore $y\notin X\setminus Z(g)$ i.e., $y\in int_XZ(g)$. Thus $y\in X_\mathscr{P}\cap int_XZ(f)\cap int_XZ(g)=int_X(X_\mathscr{P}\cap Z(f)\cap Z(g))$. Thus the set $int_X(X_\mathscr{P}\cap Z(f)\cap Z(g))$ is non-empty. It follows from Theorem \ref{Th2.6} that there exists a vertex $h\in V_\mathscr{P}(X)$, adjacent to both of $f$ and $g$. We have already obtained that $f$ and $g$ are adjacent. Hence $f-g-h-f$ is a triangle.
\end{proof}
\begin{theorem}\label{Th2.14}
	The graph $\Gamma(C_\mathscr{P}(X))\equiv\Gamma_\mathscr{P}(X)$ is triangulated if and only if for each vertex $f$, $|X_\mathscr{P}\cap int_XZ(f)|\geq 2$.
\end{theorem}
[Immediate consequence of Theorem \ref{Th2.13}.]\\\\
With the special choice $\mathscr{P}\equiv$ the ideal of all closed sets in $X$, we get the following particular case of Theorem \ref{Th2.14}.
\begin{theorem}\label{Th2.15}
	The zero-divisor graph of $C(X)$ is triangulated if and only if $X$ does not contain any isolated point.
\end{theorem}
Theorem \ref{Th2.15} was proved independently in \cite{Azarpanah} [Proposition $2.1(ii)$].\\\\
We exploit Theorem \ref{Th2.14}, to prove the following sufficient condition for $\Gamma_\mathscr{P}(X)$ to be triangulated.
\begin{theorem}
	Suppose $cl_X(X_\mathscr{P})\notin\mathscr{P}$ and each one-pointic set is a member of $\mathscr{P}$. Then $\Gamma_\mathscr{P}(X)$ is triangulated.
\end{theorem}
\begin{proof}
	It is easy to check that for any $f\in C_\mathscr{P}(X)$, $cl_X(X_\mathscr{P})\cap int_XZ(f)\subset cl_X (X_\mathscr{P}\cap int_XZ(f))$. In view of Theorem \ref{Th2.14}, it suffices to check that $X_\mathscr{P}\cap int_XZ(f)$ is an infinite set. If possible let $X_\mathscr{P}\cap int_XZ(f)=\{x_1,x_2,...,x_n\}$, a finite set.  Then it follows from the above inclusion relation that $cl_X (X_\mathscr{P})\cap int_XZ(f)=\{x_1,x_2,...,x_m\}$ for some $m\leq n$ [abusing notation]. We can write now: $cl_X(X_\mathscr{P})=cl_X(X\setminus Z(f))\cup \{x_1,x_2,...,x_m\}$. Since $f\in C_\mathscr{P}(X)$, $cl_X(X\setminus Z(f))\in\mathscr{P}$ and since each one-pointic set is a member of $\mathscr{P}$, it follows that each finite set is a member of $\mathscr{P}$. Consequently then $cl_X(X_\mathscr{P})\in\mathscr{P}$, a contradiction.
\end{proof}
We want to record the following two special cases of the last Theorem.
\begin{theorem}
	If $X$ is locally compact and non-compact, then the zero-divisor graph of $C_K(X)$ is triangulated.
\end{theorem}
\begin{theorem}
	Suppose $X$ is locally pseudocompact and non-\\pseudocompact. Then the zero-divisor graph of the ring $C_\psi(X)$ of all continuous functions with pseudocompact support is triangulated.
\end{theorem}
We shall now determine, when does $\Gamma_\mathscr{P}(X)$ become hypertriangulated. The following result is a straightforward consequence of Theorem \ref{Th2.6}.
\begin{theorem}\label{Th2.19}
	$\Gamma_\mathscr{P}(X)$ is hypertriangulated if and only if for any edge $f-g$, $int_X(X_\mathscr{P}\cap Z(f)\cap Z(g))\neq\emptyset$.
\end{theorem}
The next proposition is a sufficient condition for the hypertriangulatedness of $\Gamma_\mathscr{P}(X)$.
\begin{theorem}\label{Th2.20}
	If $cl_X(X_\mathscr{P})\notin\mathscr{P}$, then $\Gamma_\mathscr{P}(X)$ is hypertriangulated.
\end{theorem}
\begin{proof}
	It follows from Corollary \ref{cor2.4} that each non-zero element of $C_\mathscr{P}(X)$ is a vertex of $\Gamma_\mathscr{P}(X)$. Therefore if $f-g$ is an edge in this graph, then $f\neq 0$ and $g\neq 0\implies f^2+g^2\neq 0$. Hence $f^2+g^2$ is a vertex. Consequently by Theorem \ref{Th2.3}, $int_XZ(f^2+g^2)\cap X_\mathscr{P}\neq\emptyset$. This means that $int_X(Z(f)\cap Z(g)\cap X_\mathscr{P})\neq\emptyset$. Hence by Theorem \ref{Th2.19}, $\Gamma_\mathscr{P}(X)$ becomes hypertriangulated.
\end{proof}
We would like to mention in this context, the following result proved in \cite{Azarpanah}, Proposition $2.1(iii)$.
\begin{theorem}
	If $|X|>1$, then the zero-divisor graph of $C(X)$ is hypertriangulated if and only if $X$ is a connected middle $P$-space.
\end{theorem}
Since $\mathbb{R}$ is not a middle $P$-space, it follows from the last Theorem that the zero-divisor graph of $C(\mathbb{R})$ is not hypertriangulated. On putting $\mathscr{P}\equiv$ the ideal of all closed sets in $\mathbb{R}$, this reads: the zero-divisor graph of $C_\mathscr{P}(\mathbb{R})$ is not hypertriangulated. We note that with this special choice of $\mathscr{P}$, $cl_\mathbb{R}(\mathbb{R}_\mathscr{P})=\mathbb{R}\in\mathscr{P}$. Thus the condition of the Theorem \ref{Th2.20} may not hold good without the hypothesis $cl_X(X_\mathscr{P})\notin\mathscr{P}$.
\section{Cycles in $\Gamma_\mathscr{P}(X)$}\label{Sec3}
For any graph $G$, obviously the length of the smallest cycle in $G$, i.e., $gr(G)\geq 3$.
\begin{theorem}
	$3\leq gr(\Gamma_\mathscr{P}(X))\leq 4$.
\end{theorem}
\begin{proof}
	Let $f\in V_\mathscr{P}(X)$. Then there exists $g\in V_\mathscr{P}(X)$ such that $f.g=0$. Thus we always have a square in $\Gamma_\mathscr{P}(X)$: $f-g-2f-2g-f$. Therefore, $gr(\Gamma_\mathscr{P}(X))\leq 4$. 
\end{proof}
\begin{theorem}\label{Th3.2}
	If $X_\mathscr{P}$ contains atleast three points, then $gr(\Gamma_\mathscr{P}(X))=3$.
\end{theorem}
\begin{proof}
	We have to find an $f\in V_\mathscr{P}(X)$ such that $|X_\mathscr{P}\cap int_XZ(f)|\geq 2$ and the rest follows from Theorem \ref{Th2.13}. Let $x,y,z$ be distinct points in $X_\mathscr{P}$. Then by complete regularity of $X$, there exists $g\in C(X)$ such that $x\in int_XZ(g)$ and $y,z\in X\setminus Z(g)$. From Lemma \ref{Lem2.2}, there exists $f\in C_\mathscr{P}(X)$ such that $x\in X\setminus Z(f)\subset cl_X(X\setminus Z(f))\subset int_XZ(g)$. Since $y,z\in X\setminus Z(g)$, then $y,z\in X\setminus cl_X(X\setminus Z(f))=int_XZ(f)$, i.e., $y,z\in X_\mathscr{P}\cap int_XZ(f)$ and so $f\in V_\mathscr{P}(X)$.
\end{proof}
With the special choice $\mathscr{P}\equiv$ the ideal of all closed sets in $X$, we get the following particular case of Theorem \ref{Th3.2}.
\begin{theorem}\label{Th3.3}
	Whenever $X$ has atleast three points, the girth of the zero-divisor graph of $C(X)$ is $3$.
\end{theorem}
Theorem \ref{Th3.3} was proved independently in \cite{Azarpanah} [Corollary $1.3$].\\\\
The following proposition contains an exhaustive list of the length of all possible smallest cycles joining two distinct vertices.
\begin{theorem}
	Let $f,g\in V_\mathscr{P}(X)$. Then:
	\begin{enumerate}
		\item $c(f,g)=3$ if and only if $X\setminus Z(f)\cap X\setminus Z(g)=\emptyset$ and $int_X(X_\mathscr{P}\cap Z(f)\cap Z(g))\neq\emptyset$.\label{Th3.1(1)}
		\item $c(f,g)=4$ if and only if either $X\setminus Z(f)\cap X\setminus Z(g)=\emptyset$ and $int_X(X_\mathscr{P}\cap Z(f)\cap Z(g))=\emptyset$ or $X\setminus Z(f)\cap X\setminus Z(g)\neq\emptyset$ and $int_X(X_\mathscr{P}\cap Z(f)\cap Z(g))\neq\emptyset$.\label{Th3.1(2)}
		\item $c(f,g)=6$ if and only if $X\setminus Z(f)\cap X\setminus Z(g)\neq\emptyset$ and $int_X(X_\mathscr{P}\cap Z(f)\cap Z(g))=\emptyset$.
	\end{enumerate}
\end{theorem}
\begin{proof}
	\hspace{3cm}
	\begin{enumerate}
		\item Follows from Theorem \ref{Th2.6} and Theorem $2.7(\ref{Th2.7(1)})$.
		\item Let $c(f,g)=4$. If $X\setminus Z(f)\cap X\setminus Z(g)=\emptyset$, then it follows from part (\ref{Th3.1(1)}) of this Theorem that $int_X(X_\mathscr{P}\cap Z(f)\cap Z(g))=\emptyset$. On the other hand if $X\setminus Z(f)\cap X\setminus Z(g)\neq\emptyset$, then clearly $f$ and $g$ are non-adjacent. But since $c(f,g)=4$, there exists a square of the form: $f-h-g-k-f$. Thus $f$ and $g$ have a common adjacent vertex ($h$ or $k$). It follows from Theorem \ref{Th2.6} that $int_X(X_\mathscr{P}\cap Z(f)\cap Z(g))\neq\emptyset$.\\
		To prove the converse let the condition hold. Then it follows from part (\ref{Th3.1(1)}) of this Theorem that $c(f,g)>3$. Now if $X\setminus Z(f)\cap X\setminus Z(g)=\emptyset$, then $f$ and $g$ are adjacent vertices, in which case $f-g-2f-2g-f$ is a $4$-cycle containing $f$ and $g$. Hence $c(f,g)=4$ in this case. On the other hand if $X\setminus Z(f)\cap X\setminus Z(g)\neq\emptyset$ and $int_X(X_\mathscr{P}\cap Z(f)\cap Z(g))\neq\emptyset$, then this yields in view of Theorem \ref{Th2.6}, there exists a vertex $h$, adjacent to both $f$ and $g$ while $f$ and $g$ are non-adjacent. These results the $4$-cycle $f-h-g-2h-f$. Hence $c(f,g)=4$ in this case also.
		\item First assume that $c(f,g)=6$. Then it follows from (\ref{Th3.1(1)}) and (\ref{Th3.1(2)}) of this Theorem that $X\setminus Z(f)\cap X\setminus Z(g)\neq\emptyset$ and $int_X(X_\mathscr{P}\cap Z(f)\cap Z(g))=\emptyset$.\\
		Conversely let the conditions hold. Then it follows from (\ref{Th3.1(1)}) and (\ref{Th3.1(2)}) of this Theorem that $c(f,g)\neq 3$ and $c(f,g)\neq 4$, i.e., $c(f,g)>4$. Now the assumed conditions imply in view of Theorem $2.7(\ref{Th2.7(3)})$ that $d(f,g)=3$. So there exists a path $f-l-k-g$ of length $3$ joining $f$ and $g$. Surely then $f-l-k-g-2k-2l-f$ is a $6$-cycle containing $f$ and $g$. To complete the proof it remains therefore to show that there does not exist any $5$-cycle in this graph joining $f$ and $g$. We argue by contradiction. If possible let there exist a $5$-cycle which is either of the form: $f-l-k-g-h-f$, taking care of $d(f,g)=3$ or of the form: $f-l-k-h-g-f$, when there is a path of length $4$ joining $f$ and $g$. The first possibility contradicts Theorem \ref{Th2.6} while the second contradicts the observation that $d(f,g)=3$.
	\end{enumerate}
\end{proof}
The following diagrams are the graphical representations of the above Theorem.
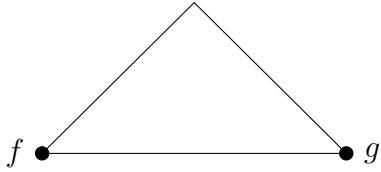
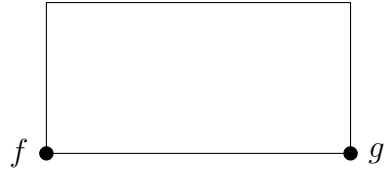
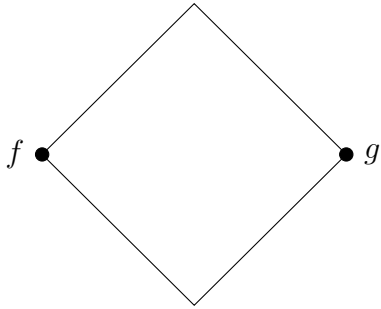
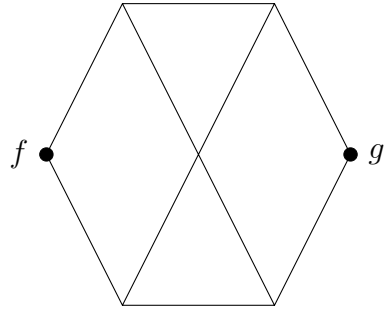
\begin{center}
	\begin{tikzpicture}
		\draw (0,0)--(2,2)--(4,0)--(0,0) (7,0)--(7,2)--(11,2)--(11,0)--(7,0);
		\coordinate[Bullet=black,label=left:$f$] (f) at (0,0);
		\coordinate[Bullet=black,label=right:$g$] (g) at (4,0);
		\coordinate[Bullet=black,label=left:$f$] (f) at (7,0);
		\coordinate[Bullet=black,label=right:$g$] (g) at (11,0);
		\node (start) at (2,-1) {$X\setminus Z(f)\cap X\setminus Z(g)=\emptyset$};
		\node (start) at (9,-1) {$X\setminus Z(f)\cap X\setminus Z(g)=\emptyset$};
		\node (start) at (2,-1.5) {$int_X(X_\mathcal{P}\cap Z(f)\cap Z(g))\neq\emptyset$};
		\node (start) at (9,-1.5) {$int_X(X_\mathcal{P}\cap Z(f)\cap Z(g))=\emptyset$};
		\draw (0,-5)--(2,-3)--(4,-5)--(2,-7)--(0,-5) (7,-5)--(8,-3)--(10,-3)--(11,-5)--(10,-7)--(8,-7)--(7,-5) (8,-3)--(10,-7) (10,-3)--(8,-7);
		\coordinate[Bullet=black,label=left:$f$] (f) at (0,-5);
		\coordinate[Bullet=black,label=right:$g$] (g) at (4,-5);
		\coordinate[Bullet=black,label=left:$f$] (f) at (7,-5);
		\coordinate[Bullet=black,label=right:$g$] (g) at (11,-5);
		\node (start) at (2,-8) {$X\setminus Z(f)\cap X\setminus Z(g)\neq\emptyset$};
		\node (start) at (9,-8) {$X\setminus Z(f)\cap X\setminus Z(g)\neq\emptyset$};
		\node (start) at (2,-8.5) {$int_X(X_\mathcal{P}\cap Z(f)\cap Z(g))\neq\emptyset$};
		\node (start) at (9,-8.5) {$int_X(X_\mathcal{P}\cap Z(f)\cap Z(g))=\emptyset$};
	\end{tikzpicture}
\end{center}
The next Corollary directly follows from these diagrams.
\begin{corollary}
	\hspace{3cm}
	\begin{enumerate}
		\item Each chord-less cycle in $\Gamma_\mathscr{P}(X)$ is of length $3$ or $4$.
		\item Every edge in $\Gamma_\mathscr{P}(X)$ is either  an edge of a triangle or an edge of a square.
	\end{enumerate}
\end{corollary}
\section{Relation between dominating number, chromatic number and clique number of $\Gamma_\mathscr{P}(X)$}\label{Sec4}
We start observing the fact which relates the clique number of $\Gamma_\mathscr{P}(X)$ and the cellularity of the space $X_\mathscr{P}$ of $X$.
\begin{theorem}
	$\omega(\Gamma_\mathscr{P}(X))=c(X_\mathscr{P})$.
\end{theorem}
\begin{proof}
	We shall first show that for an arbitrary complete subgraph $H$ of $\Gamma_\mathscr{P}(X)$, $|H|\leq c(X_\mathscr{P})$ and this will imply that $\omega(\Gamma_\mathscr{P}(X))\leq c(X_\mathscr{P})$. If $V(H)$ is the set of all vertices in $H$, then for each $f\in V(H)$, $X\setminus Z(f)$ is a non-empty open set contained in $X_\mathscr{P}$. Consequently by Theorem $2.7(\ref{Th2.7(1)})$, $\{X\setminus Z(f):f\in V(H)\}\equiv\mathscr{B}$ becomes a cellular family in $X_\mathscr{P}$. Hence $|H|=|\mathscr{B}|\leq c(X_\mathscr{P})$. To prove the reverse inequality, $c(X_\mathscr{P})\leq \omega(\Gamma_\mathscr{P}(X))$, it suffices to show for an arbitrary selected cellular family $\mathscr{B}$ in $X_\mathscr{P}$ that $|\mathscr{B}|\leq\omega(\Gamma_\mathscr{P}(X))$. Indeed for each set $B$ in the family $\mathscr{B}$, choose a point $x_B\in B$. Then by using Lemma \ref{Lem2.2} we can find out an $f_B\in C_\mathscr{P}(X)$ such that $x_B\in X\setminus Z(f_B)\subset cl_X (X\setminus Z(f_B))\subset B$. If $H$ is a subgraph of $\Gamma_\mathscr{P}(X)$ where the set of vertices is $\{f_B:B\in\mathscr{B}\}$, then the cellularity of $\mathscr{B}$ conjoined with Theorem $2.7(\ref{Th2.7(1)})$, therefore ensures that $H$ is a complete subgraph of $\Gamma_\mathscr{P}(X)$. This implies that $|\mathscr{B}|=|H|\leq\omega(\Gamma_\mathscr{P}(X))$.
\end{proof}
\begin{corollary}
	The clique number of the zero-divisor graph of $C(X)$ and the cellularity of $X$ are identical.
\end{corollary}
\begin{proof}
	This follows on choosing $\mathscr{P}\equiv$ the ideal of all closed sets in $X$.
\end{proof}
[This result is proved independently in \cite{Azarpanah}, Proposition $3.1$.]\\\\
Since the chromatic number of any graph is not less than its clique number, the following proposition is immediate:
\begin{theorem}\label{Cor4.3}
	$\chi(\Gamma_\mathscr{P}(X))\geq c(X_\mathscr{P})$.
\end{theorem}
The weight of a topological space $X$, denoted by $w(X)$, is the smallest of the cardinal numbers of the open bases for $X$.
\begin{theorem}
	$dt(\Gamma_\mathscr{P}(X))\leq w(X_\mathscr{P})$.
\end{theorem}
\begin{proof}
	Let $\mathscr{B}$ be an open base for the subspace $X_\mathscr{P}$. It suffices to find out a dominating set $D$ in $\Gamma_\mathscr{P}(X)$ with $|D|\leq|\mathscr{B}|$. Let $B\in\mathscr{B}$ such that $B\neq\emptyset$, $B\neq X_\mathscr{P}$. Fix $x_B\in B$. Then by Lemma \ref{Lem2.2}, there exists $f_B\in C_\mathscr{P}(X)$ such that $x_B\in X\setminus Z(f_B)\subset cl_X(X\setminus Z(f_B))\subset B$. Since $B\subsetneqq X_\mathscr{P}$, there exists a point $y\in X_\mathscr{P}$ such that $y\notin B$. It follows from the last inclusion relation that $y\in X-cl_X(X\setminus Z(f_B))=int_XZ(f_B)$. Consequently, $y\in X_\mathscr{P}\cap int_XZ(f_B)$. Hence from Theorem \ref{Th2.3}, $f_B\in V_\mathscr{P}(X)$. Let $D=\{f_B:B\in\mathscr{P}\}$. We claim that $D$ is a dominating set in $\Gamma_\mathscr{P}(X)$. Towards that claim choose $f\in V_\mathscr{P}(X)$. Then from Theorem \ref{Th2.3}, $X_\mathscr{P}-cl_X(X\setminus Z(f))$ is a non-empty open set in $X_\mathscr{P}$. Therefore, there exists $B\in\mathscr{B}$ such that $B\subset X_\mathscr{P}-cl_X(X\setminus Z(f))$. Consequently $f_B.f=0$. It is clear that $|D|\leq |\mathscr{B}|$, since the map: $ \left.
	\begin{array}{ll}
		\mathscr{B} & \to D \\
		B & \mapsto f_B
	\end{array}
	\right \}$ is onto $D$.
\end{proof}
\begin{definition}
	In a graph $G$, two distinct vertices $u$ and $v$ are called orthogonal if $u$ and $v$ are adjacent and there is no third vertex adjacent to $u$ and $v$ both. In this case, we write $u\perp v$. $G$ is called complemented if given a vertex $u$ in $G$, there exists a vertex $v$ in $G$ such that $u\perp v$. 
\end{definition}
The following result is a consequence of Theorem \ref{Th2.6} and Theorem \ref{Th2.7}.
\begin{theorem}\label{Th4.5}
	$\Gamma_\mathscr{P}(X)$ is complemented if and only if given $f\in V_\mathscr{P}$, there exists $g\in V_\mathscr{P}$ such that $X\setminus Z(f)\cap X\setminus Z(g)=\emptyset$ and $int_X(X_\mathscr{P}\cap Z(f)\cap Z(g))=\emptyset$.
\end{theorem}
It follows from Theorem \ref{Th2.20} that if $cl_X(X_\mathscr{P})\notin\mathscr{P}$, then $\Gamma_\mathscr{P}(X)$ is not complemented. This means that for a complemented graph $\Gamma_\mathscr{P}(X)$, $cl_X(X_\mathscr{P})\in\mathscr{P}$. On choice of $\mathscr{P}\equiv$ the ideal of all closed sets in $X$, we observe that $cl_X(X_\mathscr{P})\in\mathscr{P}$ and $C_\mathscr{P}(X)=C(X)$. Therefore the zero-divisor graph of $C(X)$ is a candidate for a complemented graph. Indeed the following fact is proved in \cite{Azarpanah}, Corollary $2.5$.
\begin{theorem}
	The zero-divisor graph of $C(X)$ is complemented if and only if the space of minimal prime ideals in $C(X)$ is compact.
\end{theorem}
We are now going to establish that this Theorem can be deduced as a special case of the more general Theorem which says that for any ideal $\mathscr{P}$ of closed sets in $X$, $\Gamma_\mathscr{P}(X)$ is complemented if and only if the space of minimal prime ideals of the ring $C_\mathscr{P}(X)$ is compact. We need a little bit of technicalities to arrive at this result.\\
We reproduce from \cite{Henriksen} the following basic information related to the space of minimal prime ideals of a commutative ring $A$ (possibly without identity), which is further reduced in the sense that there does not exist any non-zero nilpotent member of $A$. A prime ideal $P$ of $A$ is called a minimal prime ideal if there does not exist any prime ideal $Q$ of $A$ such that $Q\subsetneqq P$. It is easy to check on using Zorn's Lemma in a straight forward manner that if $P$ is a prime ideal in $A$, then there is a minimal prime ideal $Q$ in $A$ such that $Q\subset P$. Let $\mathcal{P}(A)$ be the set of all minimal prime ideals in $A$. For any subset $S$ of $A$, the hull of $S$, denoted by $h(S)$, is defined as $h(S)=\{P\in \mathcal{P}(A):S\subset P\}$. If $S$ is a single point $=\{s\}$, then we write $h(s)$ instead of $h(\{s\})$. It turns out that the family $\{h(a):a\in A\}$ is a base for the closed sets for some topology on $\mathcal{P}(A)$. $\mathcal{P}(A)$ equipped with this topology is often called the space of minimal prime ideals in $A$. For any subset $S$ of $A$, the set $\mathcal{A}(S)=\{b\in A:bS=0\}$ is called the annihilator of $S$. We just state the following results which are already proved in \cite{Henriksen}.
\begin{theorem}\label{Th4.7}
	For each member $a$ of $A$, $h(\mathcal{A}(a))=\mathcal{P}(A)\setminus h(a)$. In particular therefore $h(a)$ is a clopen set in $\mathcal{P}(A)$. Consequently $\mathcal{P}(A)$ becomes a zero-dimensional space and it is easy to prove that $\mathcal{P}(A)$ is also Hausdorff.
\end{theorem}
\begin{theorem}\label{Th4.7(A)}
	For any subset $S$ of $A$, $\mathcal{A}(S)=$ the intersection of all minimal prime ideals in $A$, which contains $S$.
\end{theorem}
\begin{theorem}\label{Th4.8}
	For any two points $x,y$ in $A$, $\mathcal{A}(\mathcal{A}(x))=\mathcal{A} (y)$ if and only if $h(x)=h(\mathcal{A}(y))$. 
\end{theorem}
$A$ is said to satisfy the annihilator condition or is called an a.c. ring if for $x,y\in A$, there exists $z\in A$ such that $\mathcal{A}(z)=\mathcal{A}(x) \cap\mathcal{A}(y)$. Given $x\in A$, an element $x'\in A$ is called a complement of $x$ if $\mathcal{A}(\mathcal{A}(x'))=\mathcal{A}(x)$. It is easy to see that if $x'$ is a complement of $x$, then $x$ is a complement of $x'$.\\
The following result relates the existence of complement of each element of $A$ with the compactness of the space of minimal prime ideals of $A$. 
\begin{theorem}
	The following statements are equivalent for the ring $A$.
	\begin{enumerate}
		\item The space $\mathcal{P}(A)$ is compact and $A$ is an a.c. ring.\label{Th4.9(1)}
		\item Each member of $A$ has a complement.
	\end{enumerate}
\end{theorem}
We reproduce the following results which appeared as Lemma $1.2$ and Theorem $0.1$ in \cite{Acharyya2} and also in \cite{MN}, Example $10$, $4.9$, page $66$.
\begin{theorem}\label{Th4.10}
	If $M$ is an ideal of $A$ such that the quotient ring $A/M$ is a field, then $M$ is a maximal ideal in $A$.
\end{theorem}
\begin{theorem}
	If $A=A^2\equiv\{\sum_{i=1}^{n}a_i.b_i: a_i,b_i\in A, n\in\mathbb{N}\}\equiv$ the internal direct product of $A$ with itself, then every maximal ideal in $A$ is prime.
\end{theorem}
With the convention $|X_\mathscr{P}|\geq 2$ made in \ref{Con2.9}, we first show that the space $\mathcal{P}(C_\mathscr{P}(X))$ of minimal prime ideals of the ring $C_\mathscr{P}(X)$ is non-empty. For that purpose, choose a point $x\in X_\mathscr{P}$.
\begin{theorem}
	The ideal $M^\mathscr{P}_x=\{f\in C_\mathscr{P}(X):f(x)=0\}$ is a maximal ideal in $C_\mathscr{P}(X)$.
\end{theorem}
\begin{proof}
	Let $t:C_\mathscr{P}(X)\to\mathbb{R}$ be the map defined by $t(f)=f(x)$. Then $t$ is a ring homomorphism. We assert that $t$ is onto $\mathbb{R}$; indeed by Lemma \ref{Lem2.2}, there exists $f\in C_\mathscr{P}(X)$ such that $f(x)\neq0$. Consequently given $r\in\mathbb{R}$, the function $\frac{r}{f(x)}.f\in C_\mathscr{P}(X)$ and $t(\frac{r}{f(x)}.f)=r$. Therefore the residue class ring of $C_\mathscr{P}(X)$ modulo the kernel of $t$ becomes isomorphic to $\mathbb{R}$. It follows that $C_\mathscr{P}(X)/ker(t)$ is a field. Hence by Theorem \ref{Th4.10}, $M^\mathscr{P}_x=ker(t)$, is a maximal ideal in $C_\mathscr{P}(X)$.
\end{proof}
Since for any $g\in C_\mathscr{P}(X)$, $g^\frac{1}{3}$ and $g^\frac{2}{3}$ also belong to $C_\mathscr{P}(X)$ and $g=g^\frac{1}{3}.g^\frac{2}{3}$, it follows that $C_\mathscr{P}(X)$ is identical with the internal direct product with itself. Hence in view of Theorem \ref{Th4.10}, we can make the following comment.
\begin{remark}
	$M^\mathscr{P}_x$ is a prime ideal in $C_\mathscr{P}(X)$. Consequently,\\ $\mathcal{P}(C_\mathscr{P}(X))\neq\emptyset$.
\end{remark}
Before proceeding further, we make the simple observation that $C_\mathscr{P}(X)$ is an a.c. ring, because for $f,g\in C_\mathscr{P}(X)$, $\mathcal{A}(f)\cap\mathcal{A}(g)=\mathcal{A}(f^2+g^2)$.\\
The following subsidiary result will be helpful to us towards proving the main result of this section.
\begin{theorem}\label{Th4.14}
	Let $f,g\in C_\mathscr{P}(X)$. Then
	\begin{enumerate}
		\item $h(g)\supset h(\mathcal{A}(f))$ if and only if $f.g=0$.
		\item $h(g)\subset h(\mathcal{A}(f))$ if and only if $X_\mathscr{P}\cap int_XZ(f)\cap int_XZ(g)=\emptyset$.
	\end{enumerate}
\end{theorem}
[Lemma $5.4$ in \cite{Henriksen} is a special case of this Theorem on choosing $\mathscr{P}\equiv$ the ideal of all closed sets in $X$.]
\begin{proof}
	\hspace{3cm}
	\begin{enumerate}
		\item If $f.g=0$, then it is clear that $g\in\mathcal{A}(f)$ and consequently $h(g)\supset h(\mathcal{A}(f))$. Conversely let $h(g)\supset h(\mathcal{A}(f))$. Then we claim in view of Theorem \ref{Th4.7(A)} that $g\in\mathcal{A}(f)$ and hence $f.g=0$.
		\item From Theorem \ref{Th4.7}, we have $h(\mathcal{A}(f))=\mathcal{P}(C_\mathscr{P}(X))\setminus h(f)$. Therefore $h(g)\subset h(\mathcal{A}(f))$ if and only if $h(g)\subset \mathcal{P}(C_\mathscr{P}(X))\setminus h(f)$, this holds if and only if $h(g)\cap h(f)=\emptyset$ and this is the case when and only when $h(f^2+g^2)=\emptyset$ meaning that $\mathcal{P}(C_\mathscr{P}(X))\setminus h(f^2+g^2)=\mathcal{P}(C_\mathscr{P}(X))$, which is the same thing in view of Theorem \ref{Th4.7} as $h(\mathcal{A}(f^2+g^2))=\mathcal{P} (C_\mathscr{P}(X))$. This means that $\mathcal{A}(f^2+g^2)\subset P$ for each $P\in\mathcal{P}(C_\mathscr{P}(X))$, equivalently $\mathcal{A}(f^2+g^2)=\{0\}$, because $C_\mathscr{P}(X)$ is a reduced ring where the intersection of all minimal prime ideals is the zero ideal. Now $\mathcal{A}(f^2+g^2)=\{0\}$ if and only if $f^2+g^2$ is not a divisor of zero in $C_\mathscr{P}(X)$. This happens in view of Theorem \ref{Th2.3}, when and only when $int_XZ(f^2+g^2)\cap X_\mathscr{P}=\emptyset$ meaning $X_\mathscr{P}\cap int_XZ(f)\cap int_XZ(g)=\emptyset$.
	\end{enumerate}
\end{proof}
A consequence of this Theorem is as follows:
\begin{theorem}\label{Th4.15}
	Given $f\in C_\mathscr{P}(X)$, a function $f'\in C_\mathscr{P}(X)$ is a complement of $f$ in this ring if and only if $f.f'=0$ and $X_\mathscr{P}\cap int_XZ(f)\cap int_XZ(f')=\emptyset$.
\end{theorem}
\begin{proof}
	$f'$ is a complement of $C_\mathscr{P}(X)$ if and only if $\mathcal{A}(\mathcal{A}(f'))=\mathcal{A}(f)$. In view of Theorem \ref{Th4.8}, this is equivalent to the statement that $h(f')=h(\mathcal{A}(f))$. If we now apply the result of Theorem \ref{Th4.14}, then we see that the last equality is equivalent to the statements that: $f.f'=0$ and $X_\mathscr{P}\cap int_XZ(f)\cap int_XZ(f')=\emptyset$.
\end{proof}
\begin{theorem}\label{Th4.16}
	The zero-divisor graph $\Gamma_\mathscr{P}(X)$ of $C_\mathscr{P}(X)$ is complemented if and only if the space $\mathcal{P}(C_\mathscr{P}(X))$ of all minimal prime ideals of $C_\mathscr{P}(X)$ is compact.
\end{theorem}
\begin{proof}
	It follows by combining Theorem \ref{Th4.5} and Theorem \ref{Th4.15} that $\Gamma_\mathscr{P}(X)$ is complemented if and only if each function $f$ in the ring $C_\mathscr{P}(X)$ has a complement. We now apply Theorem $4.11(\ref{Th4.9(1)})$ to ensure that this last statement holds if and only if $\mathcal{P}(C_\mathscr{P}(X))$ is compact.
\end{proof}
An additional observation.
\begin{remark}
	It follows from Theorem \ref{Th2.20} that if $X$ is a locally compact non-compact space (respectively a locally pseudocompact non-pseudocompact space), then the zero-divisor of $C_K(X)$ (respectively $C_\psi(X)$) is hypertriangulated and therefore not complemented. It follows from Theorem \ref{Th4.16} that the space $\mathcal{P}(C_K(X))$ of minimal prime ideals of such a $C_K(X)$ (respectively the space $\mathcal{P}(C_\psi(X))$ of minimal prime ideals of $C_\psi(X)$) is non-compact. This is an instance of how a graph theoretic result leads to a result in topology.
\end{remark}
A complemented graph $G$ is called uniquely complemented, whenever $u\perp v$ and $u\perp w$ for any three vertices $u,v,w$ in $G$, then $v$ and $w$ are adjacent to exactly the same vertices. 
\begin{theorem}\label{Th4.20}
	If $\Gamma_\mathscr{P}(X)$ is complemented, then it is uniquely complemented.
\end{theorem}
\begin{proof}
	Let $\Gamma_\mathscr{P}(X)$ be complemented and $f\perp g$, $f\perp h$ for some $f,g,h\in V_\mathscr{P}(X)$. We claim that $\mathcal{A}(g)=\mathcal{A}(h)$, where $\mathcal{A}(g)=\{l\in C_\mathscr{P}(X):l.g=0\}$. If we establish our claim, then $g$ and $h$ are adjacent to exactly the same vertices and hence $\Gamma_\mathscr{P}(X)$ is uniquely complemented. Since $f\perp g$, by Theorem \ref{Th4.5}, $X\setminus Z(f)\cap X\setminus Z(g)=\emptyset\implies X\setminus Z(g)\cap cl_X(X\setminus Z(f))=\emptyset$ and $int_X(X_\mathscr{P}\cap Z(f)\cap Z(g))=\emptyset\implies X_\mathscr{P}-(cl_X(X\setminus Z(f))\cup cl_X(X\setminus Z(g)))=\emptyset\implies X_\mathscr{P}\subset cl_X(X\setminus Z(f))\cup cl_X(X\setminus Z(g))$. Similarly $X\setminus Z(h)\cap cl_X(X\setminus Z(f))=\emptyset$ and $X_\mathscr{P}\subset cl_X(X\setminus Z(f))\cup cl_X(X\setminus Z(g))$. Now let $l\in\mathcal{A}(g)$, then $l.g=0\implies X\setminus Z(l)\cap cl_X(X\setminus Z(g))=\emptyset\implies X\setminus Z(l)\subset cl_X(X\setminus Z(f))\implies X\setminus Z(l)\cap X\setminus Z(h)=\emptyset \implies l.h=0\implies l\in \mathcal{A}(h)$, i.e., $\mathcal{A}(g)\subset \mathcal{A}(h)$. Similarly $\mathcal{A}(h)\subset \mathcal{A}(g)$ and hence $\mathcal{A}(g)=\mathcal{A}(h)$.
\end{proof}
The next result directly follows from Theorem \ref{Th4.16} and Theorem \ref{Th4.20}.
\begin{corollary}
	The zero-divisor graph $\Gamma_\mathscr{P}(X)$ of $C_\mathscr{P}(X)$ is uniquely complemented if and only if the space $\mathcal{P}(C_\mathscr{P}(X))$ of all minimal prime ideals of $C_\mathscr{P}(X)$ is compact.
\end{corollary}
\section{The zero-divisor graph of $C^\mathscr{P}_\infty(X)$}
Let $\Gamma^\mathscr{P}_\infty(X)$ stand for the zero-divisor graph of $C^\mathscr{P}_\infty(X)$ whose set of vertices is the aggregate of all non-zero zero divisors of this ring and two distinct vertices $f$ and $g$ are adjacent if and only if $f.g=0$. Let $V^\mathscr{P}_\infty(X)$ be the set of vertices of $\Gamma^\mathscr{P}_\infty(X)$. It can be realized without any difficulty that most of the results related to the zero-divisor graph $\Gamma_\mathscr{P}(X)$ of $C_\mathscr{P}(X)$ have their analogs for the zero-divisor graph $\Gamma^\mathscr{P}_\infty(X)$ of $C^\mathscr{P}_\infty(X)$. Since the proof of these later results are also parallel to the proof of their corresponding counterparts involving $\Gamma_\mathscr{P}(X)$, as obtained in Section \ref{Sec2}, \ref{Sec3} and \ref{Sec4}, we simply omit them. However we state all these parallel facts for $\Gamma^\mathscr{P}_\infty(X)$, for our convenience.
\begin{theorem}\label{Th5.1}
	For each $f\in C^\mathscr{P}_\infty(X)$, $X\setminus Z(f)\subset X_\mathscr{P}$ and an $f\in C^\mathscr{P}_\infty(X)$ is a member of $V^\mathscr{P}_\infty(X)$ if and only if $X_\mathscr{P}\cap int_XZ(f)\neq\emptyset$.
\end{theorem}
\begin{theorem}
	An $f\in V^\mathscr{P}_\infty(X)$ is a vertex of a triangle in $\Gamma^\mathscr{P}_\infty(X)$ if and only if $|X_\mathscr{P}\cap int_XZ(f)|\geq 2$.
\end{theorem}
\begin{theorem}
	If $cl_X(X_\mathscr{P})\notin\mathscr{P}$ and every finite set in $X$ is a member of $\mathscr{P}$, then $\Gamma^\mathscr{P}_\infty(X)$ is triangulated.
\end{theorem}
\begin{theorem}
	If $cl_X(X_\mathscr{P})\notin\mathscr{P}$, then $\Gamma^\mathscr{P}_\infty(X)$ is hypertriangulated.
\end{theorem}
\begin{theorem}
	Let $f,g\in V^\mathscr{P}_\infty(X)$. Then 
	$$d(f,g)=\begin{cases}
		1 & \text{ if }f.g=0\\
		2 & \text{ if }f.g\neq 0\text{ and }int_X(X_\mathscr{P}\cap Z(f)\cap Z(g))\neq\emptyset\\
		3 & \text{ if }f.g\neq 0\text{ and }int_X(X_\mathscr{P}\cap Z(f)\cap Z(g))=\emptyset
	\end{cases}$$
\end{theorem}
\begin{theorem}
	$3\leq gr(\Gamma^\mathscr{P}_\infty(X))\leq 4$.
\end{theorem}
\begin{theorem}
	For $f,g\in V^\mathscr{P}_\infty(X)$, $c(f,g)$ is either $3$ or $4$ or $6$.
\end{theorem}
\begin{theorem}
	$C^\mathscr{P}_\infty(X)$ is an a.c. ring and $\mathcal{P}(C^\mathscr{P}_\infty(X))\neq\emptyset$. Furthermore, $\Gamma^\mathscr{P}_\infty(X)$ is uniquely complemented if and only if $\mathcal{P}(C^\mathscr{P}_\infty(X))$ is a compact space. 
\end{theorem}
\begin{theorem}
	$\omega(\Gamma^\mathscr{P}_\infty(X))=\omega(\Gamma_\mathscr{P}(X))=c( X_\mathscr{P})$.
\end{theorem}
\begin{theorem}
	$dt(\Gamma^\mathscr{P}_\infty(X))\leq w(X_\mathscr{P})$.
\end{theorem}
We now state and establish two new results in this section which connect $\Gamma_\mathscr{P}(X)$ and $\Gamma^\mathscr{P}_\infty(X)$.
\begin{theorem}
	$V_\mathscr{P}(X)$ is a dominating set in the graph $\Gamma^\mathscr{P}_\infty(X)$.
\end{theorem}
\begin{proof}
	Let $f\in V^\mathscr{P}_\infty(X)$. Then from Theorem \ref{Th5.1}, we can choose a point $x\in X_\mathscr{P}\cap int_XZ(f)$. From Lemma \ref{Lem2.2}, we can have a $g\in C_\mathscr{P}(X)$ such that $x\in X\setminus Z(g)\subset cl_X(X\setminus Z(g))\subset X_\mathscr{P}\cap int_XZ(f)$. It follows that $g\in V_\mathscr{P}(X)$ and $g.f=0$, thus $g$ is adjacent to $f$. 
\end{proof}
\begin{theorem}
	$\chi(\Gamma^\mathscr{P}_\infty(X))=\chi(\Gamma_\mathscr{P}(X))$.
\end{theorem}
\begin{proof}
	Since $\Gamma_\mathscr{P}(X)$ is a subgraph of $\Gamma^\mathscr{P}_\infty(X)$, it is plain that $\chi(\Gamma^\mathscr{P}_\infty(X))\geq\chi(\Gamma_\mathscr{P}(X))$. To prove the reverse inequality let $f\in V^\mathscr{P}_\infty(X)\setminus V_\mathscr{P}(X)$. Then $X\setminus Z(f)\neq\emptyset$. Choose $x\in X\setminus Z(f)$. As $X\setminus Z(f)\subset X_\mathscr{P}$, it is clear that $x\in X_\mathscr{P}$. By Lemma \ref{Lem2.2}, there exists $g\in C_\mathscr{P}(X)$ such that $x\in X\setminus Z(g)\subset cl_X(X\setminus Z(g))\subset X\setminus Z(f)$. Let $A_f=\{g\in V_\mathscr{P}(X): X\setminus Z(g)\subset cl_X(X\setminus Z(g))\subset X\setminus Z(f)\}$. Then $A_f\neq\emptyset$ as observed above and $g\in A_f$ implies that $f$ and $g$ are non-adjacent as $f.g\neq 0$ in this case. Now there already exists a coloring of the vertices $V_\mathscr{P}(X)$ of $\Gamma_\mathscr{P}(X)$ by $\chi(\Gamma_\mathscr{P}(X))$ many colors. We want to extend this coloring to color the entire set of vertices $V^\mathscr{P}_\infty(X)$ in $\Gamma^\mathscr{P}_\infty(X)$ in a consistent manner. Indeed for any $f\in V^\mathscr{P}_\infty(X)\setminus V_\mathscr{P}(X)$, we color $f$, by the coloring of any chosen member of $A_f$. Once this assignment of colors to the members of $V^\mathscr{P}_\infty(X)$ is proved to be consistent, it will follow that the set of vertices $V^\mathscr{P}_\infty(X)$ in $\Gamma^\mathscr{P}_\infty(X)$ can be colored by the already existing colors needed to color $V_\mathscr{P}(X)$ and hence $\chi(\Gamma^\mathscr{P}_\infty(X))\leq\chi(\Gamma_\mathscr{P}(X))$. Towards the proof of the consistency of the above method of coloring, suppose $h\in V^\mathscr{P}_\infty(X)$ and $f\in V^\mathscr{P}_\infty(X)\setminus V_\mathscr{P}(X)$ have the same color as that of $g\in V_\mathscr{P}(X)$. It suffices to show that $h$ and $f$ are non-adjacent. For that purpose we need to show first that $h$ and $g$ are non-adjacent. If $h\in V_\mathscr{P}(X)$, then surely $h$ and $g$ are non-adjacent because there is already a (consistent) coloring of $V_\mathscr{P}(X)$. On the other hand if $h\in V^\mathscr{P}_\infty(X)\setminus V_\mathscr{P}(X)$, then $g\in A_g$ and hence $g$ and $h$ are non-adjacent. Therefore, $X\setminus Z(g.h)\neq\emptyset$. Also $g\in A_f$ because $f$ is colored by the color of $g$. Hence $X\setminus Z(g)\subset X\setminus Z(f)$. Consequently, $\emptyset\neq X\setminus Z(g.h)=X\setminus Z(g)\cap X\setminus Z(h)\subset X\setminus Z(f)\cap X\setminus Z(h)=X\setminus Z(f.h)$. Thus $f.h\neq 0$ and hence $f$ and $h$ are non-adjacent.
\end{proof}
\section{$\Gamma_\mathscr{P}(X)$ for a special choice of $\mathscr{P}$}\label{Sec6}
The main result in this concluding section of the present article is to prove a Banach-Stone like theorem, which tells that for appropriate choices of $\mathscr{P}$ and $\mathscr{Q}$, a graph isomorphism between $\Gamma_\mathscr{P}(X)$ and $\Gamma_\mathscr{Q}(Y)$ will lead to an isomorphism between the rings $C_\mathscr{P}(X)$ and $C_\mathscr{Q}(Y)$. Indeed we let $\mathscr{P}$ to be the ideal of all finite sets in $X$ and write $C_\mathscr{P}(X)=C_F(X)=\{f\in C(X):X\setminus Z(f)\text{ is a finite set in }X\}$. Let $\Gamma_F(X)$ denote the zero-divisor graph of $C_F(X)$ and $V_F(X)$, the set of vertices of this graph.\\
Essentially we shall show that the ring structure of $C_F(X)$ is uniquely determined by the graph structure of $\Gamma_F(X)$. We see that with this special choice of $\mathscr{P}$, $X_\mathscr{P}=K_X\equiv$ the set of all isolated points in $X$ with the Convention \ref{Con2.9}, therefore $|K_X|\geq 2$. Since for each $f\in C_F(X)$, $Z(f)$ is a clopen subset of $X$ and therefore, $X\setminus Z(f)\subset K_X$. The following special cases of Theorem \ref{Th2.3} and Theorem \ref{Th2.7} are recorded below for our immediate need.
\begin{theorem}
	\hspace{3cm}
	\begin{enumerate}
		\item An $f(\neq 0)\in C_F(X)$ is a member of $V_F(X)$ if and only if $K_X\cap Z(f)\neq\emptyset$. Thus for $f\in C_F(X)$, $\emptyset\neq X\setminus Z(f)\subsetneqq K_X$ if and only if $f\in V_F(X)$.\label{Th6.1.1}
		\item Let $f,g\in V_F(X)$. Then:
		\begin{enumerate}
			\item $d(f,g)=1$ if and only if $X\setminus Z(f)\cap X\setminus Z(g)=\emptyset$.\label{Th6.1.2.a}
			\item $d(f,g)=2$ if and only if $X\setminus Z(f)\cap X\setminus Z(g)\neq\emptyset$ and $K_X\cap Z(f)\cap Z(g)\neq\emptyset$.\label{Th6.1.2.b}
			\item $d(f,g)=3$ if and only if $X\setminus Z(f)\cap X\setminus Z(g)\neq\emptyset$ and $K_X\cap Z(f)\cap Z(g)=\emptyset$.\label{Th6.1.2.c}
		\end{enumerate}
	\end{enumerate}
\end{theorem}
The following properties determine in some cases the centers of the graph $\Gamma_F(X)$.
\begin{theorem}\label{Th6.2}
	Let $K_X$ be finite and $f\in V_F(X)$. Then $e(f)=2$ if and only if $X\setminus Z(f)$ is one-pointic set. Therefore $c(\Gamma_F(X))\equiv$ the center of $\Gamma_F(X)=\{r1_x:x\in K_X, r\in\mathbb{R}\setminus\{0\}\}$, here $1_x$ stands for the characteristic function of the set $\{x\}$, i.e., $1_x(y)=\begin{cases}
		1 & \text{if }y=x\\
		0 & \text{otherwise}
	\end{cases}$
\end{theorem}
\begin{proof}
	We prove only the first part of this Theorem, because the second part follows immediately from the first part. Assume that $X\setminus Z(f)= \{x\}$ for some $x\in X$.  Let $g\in V_F(X)$. If $x\in Z(g)$, Then $f.g=0$ so that from Theorem $6.1(\ref{Th6.1.2.a})$, $d(f,g)=1$. Suppose now that $x\notin Z(g)$. Then $x\in X\setminus Z(g)\cap X\setminus Z(f)$. On the other hand we get from Theorem $6.1(\ref{Th6.1.1})$ that $X\setminus Z(g)\subsetneqq K_X$. Hence $X\setminus Z(f)\cup X\setminus Z(g)=X\setminus Z(g)\subsetneqq K_X$. This implies that $K_X\cap Z(f)\cap Z(g)\neq\emptyset$. It follows from Theorem $6.1(\ref{Th6.1.2.b})$ that $d(f,g)=2$. Thus $e(f)=2$.\\
	To prove the converse part suppose that $X\setminus Z(f)$ contains atleast two points $x,y(x\neq y)$. We shall find out a $g\in V_F(X)$ such that $d(f,g)=3$ and that finishes the Theorem. Indeed take $K_1=(Z(f)\cap K_X) \cup\{x\}$. Then since $y\in X\setminus Z(f)\subset K_X$, it follows that $y\notin K_1$, thus $K_1\subsetneqq K_X$. It is easy to verify that $K_X=K_1\cup (X\setminus Z(f))$. Take $g=1_{K_1}\equiv$ the characteristic function of the set $K_1$. Then $g$ is a continuous function on $X$ as $K_X$ is a finite subset and hence $g\in V_F(X)$. So $g(x)\neq 0$ and $f(x)\neq 0$ imply that $X\setminus Z(g)\cap X\setminus Z(f)\neq\emptyset$. Furthermore $Z(f)\cap Z(g)\cap K_X=\emptyset$. It follows from Theorem $6.1(\ref{Th6.1.2.c})$ that $d(f,g)=3$.
\end{proof}
Earlier we mention in Corollary \ref{Cor2.4} that whenever $cl_X(X_\mathscr{P})\notin\mathscr{P}$, then every non-zero element of $C_\mathscr{P}(X)$ is a vertex of $\Gamma_\mathscr{P}(X)$. On choosing $\mathscr{P}\equiv$ the ideal of all finite sets in $X$, the converse is also true.
\begin{theorem}
	Every non-zero element of $C_F(X)$ is a vertex of $\Gamma_F(X)$ if and only if $cl_X(K_X)\notin\mathscr{P}$. In other words, every non-zero element of $C_F(X)$ is a vertex of $\Gamma_F(X)$ if and only if $K_X$ is infinite.
\end{theorem}
\begin{proof}
	If $K_X$ is infinite, so is also $cl_X(K_X)$ and hence $cl_X(K_X)\notin\mathscr{P}$. So by Corollary \ref{Cor2.4}, every non-zero element of $C_F(X)$ is a vertex of $\Gamma_F(X)$. For the converse part, it suffices to show that whenever $K_X$ is finite, then we can find a non-zero element in $C_F(X)$ which is not a vertex of $\Gamma_F(X)$. Let $K_X$ be finite. Then $K_X$ is clopen. Consider $f=1_{K_X}\in C(X)$. Then $f\neq 0$ and $X\setminus Z(f)=K_X$, which implies that, $f\in C_F(X)$. From Theorem $6.1(\ref{Th6.1.1})$, it follows that $f\notin V_F(X)$.
\end{proof}
The next proposition sharpens the result in Theorem \ref{Cor4.3}.
\begin{theorem}\label{Th6.3}
	$\chi(\Gamma_F(X))=|K_X|$
\end{theorem}
\begin{proof}
	It follows from Theorem \ref{Cor4.3} on choosing $\mathscr{P}\equiv$ the ideal of all finite sets in $X$ that $\chi(\Gamma_F(X))\geq |K_X|$. So it remains to prove only the reversed inequality. We first observe that, because of the Convention $|K_X|\geq 2$, for any point $x\in K_X$, $1_x\in V_F(X)$. This follows directly from Theorem $6.1(\ref{Th6.1.1})$. Let $H$ be the subgraph of $\Gamma_F(X)$ with its set of vertices $V(H)=\{1_x:x\in K_X\}$. Any pair of distinct vertices in $V(H)$ are adjacent because $x,y\in K_X, x\neq y\implies 1_x.1_y=0$. So there is essentially one coloring of this graph $H$, which assigns different colors to different vertices. Let us denote the color assigned to the vertex $1_x$ by notation $x$, $x\in K_X$. We now extend the coloring on $V(H)$, to a coloring on the whole graph $\Gamma_F(X)$ without breaking the consistency with the following property that, the colors $\{x:x\in K\}$, which already exist and needed to color $V(H)$ are adequate enough to color the vertices of $\Gamma_F(X)$. Once such an extension is done, the Theorem finishes thereon. For that purpose choose $f\in V_F(X)$. So $X\setminus Z(f)$ is a non-void finite set, say $X\setminus Z(f)=\{x_1, x_2,...,x_n\}$. Let us color $f$ be any $x_i~(1\leq i\leq n)$ (Surely each $x_j\in K_X$, $1\leq j\leq n$, by Theorem $6.1(\ref{Th6.1.1})$). We claim that we have already defined a consistent coloring on $V_F(X)$. So to prove this claim choose any $g\in V_F(X)$ such that $f$ and $g$ are adjacent. This means that $f.g=0$ and hence $g(x_i)=0$ because $f(x_i)\neq 0$. Thus $x_i\notin X\setminus Z(g)$ and therefore by the above mode of coloring, $g$ is colored by some element $x$ of $K_X$ with $x\neq x_i$ for each $i\in\{1,2,...n\}$. Hence the above coloring on $V_F(X)$ is consistent.
\end{proof}
\begin{remark}
	A graph is said to be locally finite if every vertex of the graph is adjacent with only finitely many vertices. Clearly, the graph $\Gamma_F(X)$ is not locally finite, because if $f\in V_F(X)$, then there exists $g\in V_F(X)$ such that $f.g=0$ and this imply that $f.(rg)=0$ for each $r\in\mathbb{R}\setminus\{0\}$. Now if we consider the topological space $X$ such that $K_X$ is finite with $|K_X|\geq 2$, then $\Gamma_F(X)$ is an example of an infinite graph which is not locally finite but $\Gamma_F(X)$ is finitely colorable.
\end{remark}
We are now ready to prove the main result of this section.
\begin{theorem}\label{Th6.4}
	Let $X$ and $Y$ be two Tychonoff spaces with $|K_X|\geq 2$ and $|K_Y|\geq 2$. Then the ring $C_F(X)$ is isomorphic to the ring $C_F(Y)$ if and only if the graph $\Gamma_F(X)$ is isomorphic to the graph $\Gamma_F(Y)$.
\end{theorem}
\begin{proof}
	If $C_F(X)$ is isomorphic to $C_F(Y)$, then it is trivial that $\Gamma_F(X)$ is isomorphic to $\Gamma_F(Y)$. Assume that there exists a graph isomorphism $\psi:\Gamma_F(X)\to\Gamma_F(Y)$ (onto $\Gamma_F(Y)$). Now if $f\in C_F(X)$, then $X\setminus Z_X(f)$ is a finite subset of $K_X$, say $X\setminus Z_X(f)=\{x_1,x_2,...,x_n\}$. So we can write $f=\sum_{i=1}^{n}f(x_i).1_{x_i}$. Formulation of an isomorphism from $C_F(X)$ onto $C_F(Y)$ is therefore guided by the following procedure: we first define a bijection on the set $K_X$ onto the set $K_Y$ and extend it linearly to a bijective map on $C_F(X)$ onto $C_F(Y)$ so that it ultimately becomes a ring isomorphism. However, we must assert before proceeding further that $|K_X|=|K_Y|$. Indeed the chromatic number of a graph is invariant under graph isomorphism, it follows that $\chi(\Gamma_F(X))=\chi(\Gamma_F(Y))$. Hence from Theorem \ref{Th6.3}, $|K_X|=|K_Y|$. Now choose $x\in K_X$, then $1_x\in V_F(X)$, as observed in the beginning of proof of Theorem \ref{Th6.3}. Consequently $\psi(1_x)\in V_F(Y)$ and let $g=\psi(1_x)$. It follows from Theorem $6.1(\ref{Th6.1.1})$ that $Y\setminus Z_Y(g)\subsetneqq K_Y$. We claim that $Y\setminus Z_Y(g)$ is a singleton.\\
	If $K_X$ is a finite set, $K_Y$ is also a finite set and from Theorem \ref{Th6.2}, we get that $e(1_x)=2$. Since $\psi$ is a graph isomorphism it follows that $e(\psi(1_x))=2$, i.e., $e(g)=2$. We apply once again Theorem \ref{Th6.2}, to ascertain that $Y\setminus Z_Y(g)$ is a singleton.\\
	Assume therefore that $K_X$ is an infinite set (and therefore $K_Y$ is an infinite set). If possible let there exist two distinct points $y_1,y_2\in Y\setminus Z_Y(g)$. Since $Y\setminus Z_Y(g)\subsetneqq K_Y$, so we can choose a point $y\in K_Y\cap Z_Y(g)$. Then $g-1_y-1_{y_1}-1_{y_2}-2.1_y-g$ ia a $5$-cycle in $\Gamma_F(Y)$. As $\psi^{-1}:\Gamma_F(Y)\to\Gamma_F(X)$ is a graph isomorphism, it follows that $1_x=\psi^{-1}(g)$ is a vertex of a $5$-cycle in $\Gamma_F(X)$, say: $1_x-f_1-f_2-f_3-f_4-1_x$. Since in this cycle $f_2$ and $1_x$ are not adjacent, it follows that $f_2.1_x\neq 0$ and hence $f_2(x)\neq 0$. By an identical reasoning $f_3(x)\neq 0$. Consequently, $f_2.f_3\neq 0$, i.e., $f_2$ and $f_3$ are not adjacent -- this contradicts the adjacency of $f_2$ and $f_3$ in the last cycle.\\
	Thus we realize that for the chosen $x\in K_X$, $Y\setminus Z_Y(\psi(1_x))$ is a one-pointic set, say $Y\setminus Z_Y(\psi(1_x))=\{y\}$, for some $y\in Y$ (eventually $y\in K_Y$, as $Y\setminus Z_Y(\psi(1_x))\subset K_Y$). This implies that $\psi(1_x)=c.1_y$ for some non-zero $c\in\mathbb{R}$. We now set $\phi(x)=y$. The map $\phi: K_X\to K_Y$ thus defined without ambiguity is certainly one-to-one and is onto $K_Y$. Indeed if $y\in K_Y$ then by following the arguments adopted above and taking care of the fact that $\psi^{-1}:\Gamma_F(Y)\to\Gamma_F(X)$ is a graph isomorphism, we can show that $\psi^{-1}(1_y)=d.1_z$ for some $z\in K_X$ and $d\neq 0$ in $\mathbb{R}$. It is easy to check that, $\phi(z)=y$.\\
	Finally define the map $\Phi:C_F(X)\to C_F(Y)$ by the following rule as contemplated earlier: if $f\in C_F(X)$, then $f=\sum_{i=1}^{n}f(x_i).1_{x_i}$, where $X\setminus Z_X(f)=\{x_1,x_2,...,x_n\}$. We set $\Phi(f)=\sum_{i=1}^{n}f(x_i).1_{\phi(x_i)}$, if $f\neq 0$; and $\Phi(0)=0$. It can be proved by some routine calculations that for $f,g\in C_F(X)$, $\Phi(f+g)=\Phi(f)+\Phi(g)$ and $\Phi(f.g)=\Phi(f). \Phi(g)$. Thus $\Phi$ is a ring homomorphism which is one-to-one because for $f=\sum_{i=1}^{n}f(x_i).1_{x_i}\in C_F(X)$, $\Phi(f)=0$ implies that $f(x_i)=0$ for each $i=1,2,...,n$ and hence $f=0$. Finally for $h=\sum_{i=1}^{n}h(y_i).1_{y_i}\in C_F(Y)$ where $Y\setminus Z_Y(h)=\{y_1.y_2,...,y_n\}$. We see that $f=\sum_{i=1}^{n}h(y_i).1_{\phi^{-1}(y_i)}\in C_F(X)$ and $\Phi(f)=h$. Thus $\Phi$ is an isomorphism from $C_F(X)$ onto $C_F(Y)$.
\end{proof}
Note that the conditions $|K_X|\geq 2$ and $|K_Y|\geq 2$ can not be dropped in the statement of the above Theorem. Consider the following example:
\begin{example}
	Consider $X=[0,1]$ and $Y=[0,1]\cup\{2\}$. Then $X$ has no isolated point whereas $2$ is the only isolated point of $Y$. So $C_F(X)=\{0\}$ and $C_F(Y)=\{r.1_2: r\in\mathbb{R}\}$ and hence $C_F(X)$ and $C_F(Y)$ are not isomorphic but the zero divisor graph of both the rings $C_F(X)$ and $C_F(Y)$ are empty and hence isomorphic.
\end{example}
\bibliographystyle{plain}

}
\end{document}